\documentclass[letterpaper, 10 pt, conference]{ieeeconf}
\IEEEoverridecommandlockouts                              
\usepackage{cite}
\usepackage{amsmath,amssymb,amsfonts}
\usepackage{algorithmic}
\usepackage{graphicx}
\usepackage{textcomp}

\usepackage{amssymb,latexsym,color,amsmath,pifont,graphicx, cite}
\usepackage{subfigure}  
\usepackage{float}
\usepackage[colorlinks=true, allcolors=blue]{hyperref}
\usepackage[nameinlink]{cleveref}
\usepackage{lipsum}

\makeatletter
\newcommand*\rel@kern[1]{\kern#1\dimexpr\macc@kerna}
\newcommand*\widebar[1]{%
	\begingroup
	\def\mathaccent##1##2{%
		\rel@kern{0.8}%
		\overline{\rel@kern{-0.8}\macc@nucleus\rel@kern{0.2}}%
		\rel@kern{-0.2}%
	}%
	\macc@depth\@ne
	\let\math@bgroup\@empty \let\math@egroup\macc@set@skewchar
	\mathsurround\z@ \frozen@everymath{\mathgroup\macc@group\relax}%
	\macc@set@skewchar\relax
	\let\mathaccentV\macc@nested@a
	\macc@nested@a\relax111{#1}%
	\endgroup
}

\DeclareMathOperator*{\argmin}{argmin}

\DeclareMathOperator*{\minimize}{minimize}

\DeclareMathOperator*{\subjectto}{subject~to}

\newcommand{\reals}[1]{\mathbb{R}^{#1}}
\newcommand{\cN}{\mathcal{N}}

\DeclareMathOperator{\prox}{\mathbf{prox}}

\DeclareMathOperator{\sign}{sign}

\newcommand{\norm}[1]{\| #1 \|}    
\newcommand{\inner}[2]{\langle #1,#2\rangle}    

\newcommand{\set}[1]{\{ #1 \}}
\newcommand{\enma}[1]{\ensuremath{#1}}

\newcommand{\non}{\nonumber}

\newcommand{\asp}[1]{\!\!\! #1 \!\!\!}

\newcommand{\beq}{\begin{equation}}
\newcommand{\eeq}{\end{equation}}
\newcommand{\ba}{\begin{array}}
\newcommand{\ea}{\end{array}}
\newcommand{\bseq}{\begin{subequations}}
\newcommand{\eseq}{\end{subequations}}

\newcommand{\DefinedAs}[0]{\mathrel{\mathop:}=}

\newcommand{\rme}{\mathrm{e}}

\newcommand{\matbegin}{\left[}
\newcommand{\matend}{\right]}

\newcommand{\tbo}[2]{
		\matbegin \begin{array}{c}
				#1 \\ #2
			\end{array} \matend }

\newcommand{\zz}{z}

\newcommand{\xs}{x^\star}

\newcommand{\Tb}{\widebar{T}}
\newcommand{\xd}{\dot{x}}
\newcommand{\zd}{\dot{\zz}}

\newcommand{\cLm}{\enma{\mathcal L}_\mu}
\newcommand{\cGm}{\enma{\mathcal G}_\mu}

\newcommand{\cG}{\enma{\mathcal G}}

\newcommand{\cM}{\enma{\mathcal M}}

\newcommand{\kk}{k}

\newtheorem{theorem}{Theorem}

\newtheorem{lemma}{Lemma}
\newtheorem{remark}{Remark}

\newtheorem{corollary}{Corollary}

\newcommand\nvs{0.2cm}
\newcommand{\vsp}{\vspace*{0.01cm}}

\begin{document}
%
\title{\bf From Exponential to Finite/Fixed-Time Stability: \\ Applications to Optimization  }
\author{Ibrahim~K.~Ozaslan,~\IEEEmembership{Graduate~Student~Member,~IEEE,}
	and Mihailo~R.~Jovanovi\'{c},~\IEEEmembership{Fellow,~IEEE}
	\thanks{I.\ K.\ Ozaslan and M.\ R.\ Jovanovi\'{c} are with the Ming Hsieh Department of Electrical and Computer Engineering, University of Southern California, Los Angeles, CA 90089.
		{E-mails: ozaslan@usc.edu, mihailo@usc.edu.}
	}
}

\maketitle

\begin{abstract}
The development of finite/fixed-time stable optimization algorithms typically involves study of specific problem instances. The lack of a unified framework hinders understanding of more sophisticated algorithms, e.g., primal-dual gradient flow dynamics. The purpose of this paper is to address the following question: Given an exponentially stable optimization algorithm, can it be modified to obtain a finite/fixed-time stable algorithm? We provide an affirmative answer, demonstrate how the solution can be computed on a finite-time interval via a simple scaling of the right-hand-side of the original dynamics, and certify the desired properties of the modified algorithm using the Lyapunov function that proves exponential stability of the original system. Finally, we examine nonsmooth composite optimization problems and smooth problems with linear constraints to demonstrate the merits of our approach. 
\end{abstract}


\vspace*{0.15cm}


\section{Introduction}\label{sec.intro}
A growing body of literature views optimization algorithms as continuous-time dynamical systems. This perspective dates back to the seminal paper~\cite{arrhuruza58} and it leverages tools from control theory to provide new algorithmic insights and schemes from numerical analysis to obtain different discrete-time strategies. Examples range from studies of momentum-based accelerated techniques~\cite{suboycan16, wibwiljor16, muejor21,frarobvid23}, to analysis of primal-dual methods~\cite{quli18,dhikhojovTAC19,ozajovCDC23,holles21} and nonconvex problems~\cite{dujinleejorsinpoc17,murswekar19, oshheawu23}, and the design of optimization algorithms~\cite{cor06,muejor22,allcor23}. In particular, two classes of (scaled or normalized) gradient flow dynamics for minimizing a continuously differentiable convex function \mbox{$f$: $\reals{n}\to\reals{}$} were proposed in~\cite{cor06},
\begin{subequations}\label{eq.gf}
	\begin{IEEEeqnarray}{rcl}
	\label{eq.gf1}
	\xd  
	& ~\in~ &
	-\nabla f(x)/\norm{\nabla f(x)}
	\\[0.15cm]
	 \xd  
	 & ~\in~ &
	 \label{eq.gf2}
	  -\sign (\nabla f(x)).
	  \end{IEEEeqnarray}
\end{subequations}
The systems in~\eqref{eq.gf} are expressed as differential inclusions since the right-hand-sides are discontinuous at the equilibrium points. For a strongly convex $f$, solutions to both of these systems converge to the minimizer of $f$ in \textit{finite-time}. Extensions to projected gradient flows are also possible~\cite{cheren20}. Besides their superb convergence in convex settings, scaled gradient flows may enjoy additional advantages in nonconvex scenarios. For example, while naive gradient flow dynamics escapes saddle points in a time that grows exponentially with the problem dimension $n$~\cite{dujinleejorsinpoc17}, the normalized gradient flow does this in finite time that is proportional to the local condition number of the objective function $f$~\cite{murswekar19}. 

The study of dynamical systems with discontinuous right-hand-sides hinders the use of standard techniques for their analysis. The following modification of the scaling factor $1/\norm{\nabla f(x)}$ makes the right-hand-side in~\eqref{eq.gf1} continuous and alleviates some challenges in stability analysis~\cite{romben20},
\beq\non
\xd  \,=\, -\nabla f(x)/\norm{f(x)}^\lambda, \qquad \lambda \,\in\,(0,1).
\eeq

\vsp

While the above gradient flows compute the optimal solution on a finite time interval, the resulting convergence times depend on the initial condition and they can grow without an upper bound. In contrast, \textit{fixed-time} stability requires a uniform upper bound on the convergence time across all initial conditions. Building on a Lyapunov-based characterization of fixed-time stability in~\cite{pol12}, the scaling factor $\norm{\nabla f(x)}^{-\lambda}$ was modified in~\cite{garpan21} to achieve fixed-time stability for various types of gradient flows used in smooth unconstrained and smooth convex optimization with equality constraints. Recent extensions also include proximal gradient flows for minimizing composite objective functions which can be expressed as the sum of a smooth and a nonsmooth term~\cite{garbargupben22, julihefen23}.

\vsp

The analysis of finite/fixed-time stability of optimization algorithms has primarily been done for specific problem instances. Results are limited to relatively simple gradient flow dynamics and have yet to be extended to more sophisticated algorithms such as the primal-dual gradient flows resulting from the proximal augmented Lagrangian of nonsmooth composite problems~\cite{dhikhojovTAC19, ozajovCDC22} or gradient flows based on the generalized Lagrangian of constrained problems~\cite{quli18, tanquli20}. 

\vsp

Existing results exploit a Polyak-Lojasiewicz (PL) type condition and uniqueness of the minimizer for achieving finite/fixed-time stability in smooth unconstrained problems, strong convexity along with certain rank assumptions in smooth convex problems with equality constraints, and strict/strong monotonicity in convex-concave cases. In this paper, we recognize a common theme between the sufficient conditions for achieving finite/fixed-time and exponential stability and address the following fundamental question: 
	\begin{itemize}
	\item \textit{Can an exponentially stable optimization algorithm be modified to obtain a finite/fixed-time stable algorithm?} 
	\end{itemize}
  

\vsp

We provide an affirmative answer to this question and demonstrate that standard exponential stability results represent an essential building block for establishing finite/fixed-time stability for a broad class of optimization algorithms. Specifically, we show how exponentially stable gradient flow dynamics with Lipschitz continuous right-hand-side can be modified to obtain finite/fixed-time stable algorithm. The modification amounts to a simple scaling of the right-hand-side of the original exponentially stable system. Owing to the modularity of our approach, this form of stability for a broad class of optimization algorithms follows from standard Lyapunov-based exponential stability results. We demonstrate the merits and the utility of our approach by establishing fixed-time stability of the scaled gradient flow dynamics associated with the aforementioned primal-dual optimization algorithms.


\vsp

We also note that study of finite/fixed-time stability in control theory is not limited to continuous-time optimization algorithms. The applications include, but are not limited to, switched systems~\cite{dulinli09}, stochastic dynamics~\cite{yuyuliyan19},  circuit design~\cite{juyuayanshi24}, optimal control~\cite{hadlee22}, and game theory~\cite{povkrsbas22}. Moreover,  finite/fixed-time stability can be achieved in discrete-time nonlinear systems either directly by design~\cite{hadlee20,adlatt20} or via the use of implicit discretization schemes~\cite{bropol15}.

\vsp

The rest of the paper is organized as follows. In Section~\ref{sec.back}, we setup the problem and provide background material. In Section~\ref{sec.main}, we present our main findings and, in Sections~\ref{sec.app}, we obtain fixed-time stable primal-dual gradient flow dynamics for two types of Lagrangians. In Section~\ref{sec.num}, we provide examples and, in Section~\ref{sec.conc}, we conclude the~paper.

	\vspace*{-1ex}
\section{Problem Formulation and background}\label{sec.back}

We consider time-invariant dynamical systems,
\beq\label{eq.dyn}
\xd  \,=\,  F(x)
\eeq
with an equilibrium point at the origin, $F(0)=0$. Additional assumptions on the dynamical generator $F$: $\reals{n}\to\reals{n}$ are introduced in Section~\ref{sec.main}. 

	\vsp
	
In what follows, we provide definitions of different stability notions along with the associated Lyapunov-based sufficient conditions.




	\vspace*{-1ex}
\subsection{Exponential stability}
The origin is the globally exponentially stable equilibrium point of~\eqref{eq.dyn} if there are positive constants $M$ and $\rho$ that are independent of $x(0)$ such that the solution $x (t)$ to~\eqref{eq.dyn} satisfies,
\beq\label{eq.exp}
\norm{x(t)}  \,\leq\, M\rme^{-\rho t}\norm{x(0)},\quad \forall \, t \,\geq\, 0.
\eeq

Lemma~\ref{lemma.exp} characterizes a necessary and sufficient condition for the exponential stability.

\vsp

\begin{lemma}[Thm. 3.11 \cite{had08}]\label{lemma.exp}
	The origin is the globally exponentially stable equilibrium point of~\eqref{eq.dyn} if and only if\footnote{While~\cite[Thm. 3.11]{had08} requires $F$ to be continuously differentiable, this condition can be relaxed to locally Lipschitz continuous vector fields; see~\cite[Sect. 3.5]{had08} for details.} there exists a continuously differentiable Lyapunov function $V$ that satisfies,
	\bseq\label{eq.lemma.ges}
	\begin{IEEEeqnarray}{rcl}
		\kk_1\norm{x}^2  
		&\;\leq\;&  
		V(x)   \,\leq\, \kk_2\norm{x}^2
		\label{eq.lemma.ges.1}
		\\[\nvs]
		\dot{V}(t)
		&\;\leq\;& 
		-\kk_3\norm{x(t)}^2,\quad \forall \, t \,\geq\, 0
		\label{eq.lemma.ges.2}
	\end{IEEEeqnarray}
	\eseq
	along the solutions of~\eqref{eq.dyn}. 
\end{lemma}  

	\vspace*{-1ex}
\subsection{Finite-time stability}
The origin is the globally finite-time stable equilibrium point of~\eqref{eq.dyn} if 
\begin{enumerate}
\item[(i)] it is stable in the sense of Lyapunov~\cite[Def.~3.1]{had08};
\item[(ii)] there is a settling-time function $T$: $\reals{n}\to(0,\infty)$ such that \mbox{$\lim_{t \, \to \, T(x(0))} x(t) = 0$}; see \cite[Def.\ 2.2]{bhaber00}. 
\end{enumerate} 

Lemma~\ref{lemma.fin} provides a Lyapunov-based characterization of finite-time stability.
\vsp

\begin{lemma}[Thm.\ 4.2,~\cite{bhaber00}]\label{lemma.fin}
	Let $V$ be a globally positive definite radially unbounded Lyapunov function with $V(0) = 0$ and let the inequality
	\beq\label{eq.lemma.fin}
	\dot{V}(t)  \,\leq\, -c V^\alpha(x(t)),\quad \forall \, t \,\geq\, 0
	\eeq
	hold along the solutions of~\eqref{eq.dyn} for some $\alpha \in (0,1)$ and $c>0$. The origin is globally finite-time stable equilibrium point of~\eqref{eq.dyn} and the settling-time is upper bounded by
	\beq
	T(x(0)) \,\leq\, \dfrac{V^{1 \,-\, \alpha}(x(0))}{c(1 \,-\, \alpha)}.
	\eeq
\end{lemma} 

	\vspace*{-1ex}
\subsection{Fixed-time stability}
The origin is the globally fixed-time stable equilibrium point of~\eqref{eq.dyn}  if 
\begin{itemize}
\item[(i)] it is globally finite-time stable; 
\item[(ii)] there exists a uniform upper bound on the settling-time for all initial conditions. 
\end{itemize}

Lemma~\ref{lemma.fxd} provides a Lyapunov-based characterization of the fixed-time stability.
\vsp

\begin{lemma}[Lemma 1,~\cite{pol12}]\label{lemma.fxd}
	Let $V$ be a globally positive definite radially unbounded Lyapunov function with $V(0) = 0$ and let the inequality
	\beq\label{eq.lemma.fxd}
	\dot{V}(t)  \,\leq\, -c_1 V^{\alpha_1}(x(t)) \,-\, c_2 V^{\alpha_2}(x(t)), \quad \forall \, t \,\geq\, 0
	\eeq
	hold along the solutions of \eqref{eq.dyn} for some  $\alpha_1\in(0,1)$ and $\alpha_2>1$. The origin is the globally fixed-time stable equilibrium point of~\eqref{eq.dyn} and, for all initial conditions, the settling time is upper bounded by 
	\beq\label{eq.lemma.fxd.time}
	\Tb \,\DefinedAs\,  \sup_{x \, \in \, \reals{n}}T(x)   \,<\,\dfrac{1}{c_1(1  \,-\, \alpha_1)} \,+\, \dfrac{1}{c_2(\alpha_2 \,-\, 1)}.
	\eeq
\end{lemma}

\section{Main results}\label{sec.main}

We now demonstrate how a globally exponentially stable system can be modified to guarantee finite/fixed-time stability of the equilibrium point. This is accomplished by a simple state-dependent scaling of the generator of the original dynamics. We use a Lyapunov-based approach to establish our results. 


We first show that if the origin is globally exponentially stable equilibrium point of the original system~\eqref{eq.dyn}, with globally Lipschitz continuous dynamical generator $F$, 
	\beq\label{eq.lip}
\norm{F(x) \,-\, F(y)}  \,\leq\, L\norm{x \,-\, y}, \quad \forall x,y  \,\in\, \reals{n}
	\eeq
then it is globally finite-time stable equilibrium point of 
\beq\label{eq.nor_dyn}
\xd  \,=\,  \sigma(x)F (x)
\eeq
where $\sigma$: $\reals{n}\to\reals{}_+\!\cup\set{0}$ is the scaling factor given by
\beq\label{eq.theorem.fin_dyn.sig}
	\sigma(x)=\left\{\ba{ll}
	0 & F(x) \, = \, 0
	\\[\nvs]
	\eta \norm{F(x)}^{-\lambda} & \text{otherwise}
	\ea\right.
\eeq
with $\eta>0$ and $\lambda\in(0,1)$. Furthermore, we prove global fixed-time stability under two additional conditions: in addition to being Lipschitz continuous, the vector field $F$ in~\eqref{eq.dyn} satisfies
\beq\label{eq.grad_lower}
\norm{F(x)} \,\geq\,  m\norm{x}^\beta, \quad \forall \, x\in\reals{n}
\eeq
for some positive parameters $m$ and $\beta$, and the scaling factor is modified to
	\beq\label{eq.theorem.fxd_dyn.sig}
	\hspace{-0.2cm}\sigma(x) \,=\, \left\{\ba{ll}
	0 & F(x) \, = \, 0
	\\[\nvs]
	\eta_1 \norm{F(x)}^{-\lambda_1} \,+\, \eta_2 \norm{F(x)}^{\lambda_2} & \text{otherwise}
	\ea\right.
	\eeq
where $\eta_1$, $\eta_2$, and $\lambda_2$ are positive parameters and $\lambda_1\in(0,1)$. When these hold, a uniform upper bound on the settling-time of system~\eqref{eq.nor_dyn} for all initial conditions can be obtained.

\vsp

In Theorems~\ref{theorem.fin}~and~\ref{theorem.fxd}, we utilize the existence of a Lyapunov function that certifies global exponential stability of system~\eqref{eq.dyn} to establish finite- and fixed-time stability of~\eqref{eq.nor_dyn} for scaling factors~\eqref{eq.theorem.fin_dyn.sig} and~\eqref{eq.theorem.fxd_dyn.sig}, respectively. Then, in Lemma~\ref{lemma.exst_uni}, we show the existence and uniqueness of solutions to~\eqref{eq.nor_dyn}.

\vsp

\begin{theorem}\label{theorem.fin}
	Let the origin be globally exponentially stable equilibrium point of~\eqref{eq.dyn} and let $(k_1,k_2,k_3)$ be the parameters of the Lyapunov function in~\eqref{eq.lemma.ges}. If $F$: $\reals{n}\to\reals{n}$ is globally Lipschitz continuous with modulus $L$, then the origin is globally finite-time stable equilibrium point of~\eqref{eq.nor_dyn} with the scaling factor~\eqref{eq.theorem.fin_dyn.sig} for any $\eta>0$ and $\lambda\in(0,1)$. Moreover,  the settling-time is upper bounded by
	\beq
	T(x(0))   \,\leq\, \dfrac{2k_2 L^\lambda}{k_3\eta\lambda}\norm{x(0)}^\lambda.
	\eeq
\end{theorem}

\vsp

\begin{proof}
	We show that the Lyapunov function $V$ that certifies the global exponential stability of~\eqref{eq.dyn} also satisfies~\eqref{eq.lemma.fin} along the solutions of~\eqref{eq.nor_dyn} with the scaling factor~\eqref{eq.theorem.fin_dyn.sig},
	\beq\non
	\ba{rcl}
	\dot{V} 
	&\asp{=}& 
	\inner{\nabla_x V(x)}{\xd}  
	\,=\, 
	\sigma(x)\inner{\nabla_x V(x)}{F(x)}
	\\[\nvs]
	&\asp{\leq}& -k_3\sigma(x)\norm{x}^2  
	\,=\,  
	-\dfrac{k_3\eta\norm{x}^2}{\norm{F(x)}^{\lambda}}
	\\[\nvs]
	&\asp{\leq}&
	-\dfrac{k_3\eta}{L^\lambda} \, \norm{x}^{2-\lambda}
	\\[\nvs]
	&\asp{\leq}&
	-\dfrac{k_3\eta}{L^\lambda k_2^{1 \,-\, \lambda/2}} \, V^{1 \,-\, \lambda/2}.
	\ea
	\eeq 
	The first inequality is obtained using~\eqref{eq.lemma.ges.2}, the second inequality follows from the Lipschitz continuity~\eqref{eq.lip}, and the last inequality is obtained using the quadratic upper bound~\eqref{eq.lemma.ges.1} on the Lyapunov function. Setting
	\beq
	\alpha  \,=\,  1 - \lambda/2, \quad c \,=\, k_3\eta/(L^\lambda k_2^\alpha)
	\eeq
	yields $\dot{V} \leq -cV^\alpha$ with $c>0$ and $\alpha\in(1/2,1)$. Finite-time stability of the origin follows from Lemma~\ref{lemma.fin}, with the following upper bound on the settling-time
	\beq\non
	T(x(0))   \,\leq\, \dfrac{2k_2^\alpha L^\lambda V^{\lambda/2}(x(0))}{k_3\eta\lambda} \,\leq\, \dfrac{2k_2L^\lambda \norm{x(0)}^\lambda}{k_3\eta\lambda}
	\eeq
	where the second inequality is obtained using~\eqref{eq.lemma.ges.1}. 
\end{proof}

%

\vsp

\begin{theorem}\label{theorem.fxd}
	Let the origin be globally exponentially stable equilibrium point of~\eqref{eq.dyn} and let $(k_1,k_2,k_3)$ be the parameters of the Lyapunov function in~\eqref{eq.lemma.ges}. If $F$: $\reals{n}\to\reals{n}$ is globally Lipschitz continuous with modulus $L$ and if it satisfies~\eqref{eq.grad_lower} with positive parameters $m$ and $\beta$, then the origin is globally fixed-time stable equilibrium point of~\eqref{eq.nor_dyn} with the scaling factor~\eqref{eq.theorem.fxd_dyn.sig} for any $\eta_1$, $\eta_2$, $\lambda_2 > 0$ and $\lambda_1\in(0,1)$. Moreover, the settling-time is uniformly upper bounded by
	\beq
	\Tb \,\leq\, \dfrac{2k_2^{1 \,+\, \beta\lambda_2/2}}{k_3}\left(\dfrac{2L^{\lambda_1}}{\eta_1\lambda_1} \,+\, \dfrac{1}{\eta_2m^{\lambda_2}\beta\lambda_2}\right).
	\eeq
\end{theorem}

\vsp

\begin{proof}
	We show that if $F$ in~\eqref{eq.dyn} satisfies~\eqref{eq.grad_lower}, then the Lyapunov function $V$ that verifies exponential stability of~\eqref{eq.dyn} also satisfies~\eqref{eq.lemma.fxd} along the solutions of system~\eqref{eq.nor_dyn} with the scaling factor~\eqref{eq.theorem.fxd_dyn.sig},
	\beq\non
	\ba{rcl}
	\dot{V} 
	&\asp{=}& 
	\inner{\nabla_x V(x)}{\xd}  
	\,=\, 
	\sigma(x)\inner{\nabla_x V(x)}{F(x)}
	\\[\nvs]
	&\asp{\leq}& -k_3\sigma(x)\norm{x}^2  
	\\[\nvs]
	&\asp{=}&
	-k_3\left(\eta_1\dfrac{\norm{x}^2}{\norm{F(x)}^{\lambda_1}} \,+\, \eta_2\norm{x}^2\norm{F(x)}^{\lambda_2}\right)
	\\[0.4cm]
	&\asp{\leq}&
	-\dfrac{k_3\eta_1}{L^{\lambda_1}}\norm{x}^{2 \,-\,  \lambda_1} \,-\,k_3 \eta_2m^{\lambda_2}\norm{x}^{2+\beta\lambda_2}
	\\[0.4cm]
	&\asp{\leq}&
	-\dfrac{k_3\eta_1}{L^{\lambda_1}k_2^{1 \,-\, \lambda_1/2}}V^{1 \,-\,  \lambda_1/2} \,-\,\dfrac{k_3 \eta_2m^{\lambda_2}}{k_2^{1 \,+\, \beta\lambda_2/2}}V^{1+\beta\lambda_2/2}.
	\ea
	\eeq 
	The first inequality is obtained using~\eqref{eq.lemma.ges.2}, the first term in the second inequality is upper bounded using the Lipschitz continuity~\eqref{eq.lip}, the second term is upper bounded using~\eqref{eq.grad_lower}, and the last inequality is obtained using the quadratic upper bound~\eqref{eq.lemma.ges.1}. Setting
	\beq\label{eq.proof.fxd.par}
	\ba{rclcrcl}
	\alpha_1  &\asp{=}& 1  \,-\, \lambda_1/2 && c_1  &\asp{=}& k_3\eta_1/(L^\lambda k_2^{\alpha_1})
	\\[\nvs]
	\alpha_2  &\asp{=}& 1 \,+\, \beta\lambda_2/2, && c_2  &\asp{=}& k_3 \eta_2m^{\lambda_2}/k_2^{\alpha_2} 
	\ea
	\eeq 
	yields $\dot{V} \leq -c_1V^{\alpha_1} - c_2V^{\alpha_2}$ with $c_1,c_2 >0$, $\alpha_1\in(1/2,1)$, and $\alpha_2>1$. Thus, Lemma~\ref{lemma.fxd} implies fixed-time stability of the origin and substituting~\eqref{eq.proof.fxd.par} into~\eqref{eq.lemma.fxd.time} results in the following uniform upper bound on the settling-time
	\beq\non
	\Tb \,\leq\, \dfrac{2L^{\lambda_1} k_2^{\alpha_1}}{k_3\eta_1\lambda_1} \,+\, \dfrac{2k_2^{\alpha_2}}{k_3 \eta_2m^{\lambda_2}\beta\lambda_2}.
	\eeq
	Since $\alpha_1 \leq \alpha_2$, replacing $\alpha_1$ by $\alpha_2$ completes the proof.
\end{proof}

\vsp

\begin{remark}
	In Theorems~\ref{theorem.fin} and~\ref{theorem.fxd}, parameters $(k_1,k_2,k_3)$ of the Lyapunov function are not needed for the algorithm design but only to derive an upper bound on the settling-time.
\end{remark} 

\vsp

	\begin{remark}
Having a Lipschitz continuous Lyapunov function that certifies the exponential stability of~\eqref{eq.dyn} suffices\footnote{Let $\nabla_xV$ in~\eqref{eq.lemma.ges.2} be Lipschitz continuous function with modulus $\ell$. Then, Cauchy-Schwartz inequality applied to~\eqref{eq.lemma.ges.2} yields 
	$$
	k_3\norm{x}^2   \,\leq\,  |\dot{V}| \,\leq\, \norm{\nabla_x V(x)}\norm{F(x)}  \,\leq\,\ell\norm{x}\norm{F(x)}
	$$
i.e., the vector field $F$ satisfies~\eqref{eq.grad_lower} with $m = k_3/\ell$ and $\beta = 1$.} to satisfy the additional condition~\eqref{eq.grad_lower} on the dynamic generator $F$ in Theorem~\ref{theorem.fxd}. This sufficient condition can be readily verified in numerous cases; see Section~\ref{sec.app}.
	\end{remark} 
	
\vsp

	\begin{remark}
Existing finite/fixed-time stability results for unconstrained minimization problems~\cite{romben20,garpan21,garbargupben22,julihefen23} assume that the (proximal) gradient flow, i.e., when $F(x) = -\nabla f(x)$ (or $F(x) = -(x - \prox_{\mu g}(x - \mu\nabla f(x)))$), has a unique equilibrium point and that the objective function $f(x)$ (or $f(x)+ g(x)$) satisfies the (proximal) PL condition. 	The equivalence between the (proximal) PL and the (proximal) error bound conditions~\cite{karnutsch16} implies that under these two assumptions the dynamical generator $F$ satisfies~\eqref{eq.grad_lower}.
	\end{remark} 

	\vsp

In Lemma~\ref{lemma.exst_uni}, we establish the existence and uniqueness of solutions to system~\eqref{eq.nor_dyn} with scaling factor~\eqref{eq.theorem.fin_dyn.sig} or~\eqref{eq.theorem.fxd_dyn.sig}. Lemma~\ref{lemma.exst_uni} relies on~\cite[Prop.~1]{garbargupben22} which shows that a solution to~\eqref{eq.nor_dyn} is also a solution to~\eqref{eq.dyn} under a suitable time parameterization, i.e., that the two dynamical systems are topologically equivalent~\cite[Sect.\ 3.1]{per13}; see~\cite[Prop.~7]{murswekar19} for a similar result. 

\vsp

\begin{lemma}\label{lemma.exst_uni}
	Let the origin be globally exponentially stable equilibrium point of system~\eqref{eq.dyn} with a globally Lipschitz continuous $F$. 
	If the scaling factor $\sigma$ is given by either~\eqref{eq.theorem.fin_dyn.sig} or~\eqref{eq.theorem.fxd_dyn.sig}, then there exists a unique solution to the modified system~\eqref{eq.nor_dyn} for any initial condition $x(0)$.
\end{lemma}

\vsp

\begin{proof}
	The proof is a straightforward application of~\cite[Prop. 1]{garbargupben22}. Since the scaling factors~\eqref{eq.theorem.fin_dyn.sig} and~\eqref{eq.theorem.fxd_dyn.sig} are continuous in~$\reals{n}\setminus\set{0}$ and since $F$ is Lipschitz continuous and the exponents in~\eqref{eq.theorem.fin_dyn.sig} and~\eqref{eq.theorem.fxd_dyn.sig} satisfy $\lambda,\lambda_1\in(0,1)$, we have $\lim_{x\to0}\sigma(x)F(x) = 0$. Thus, the right-hand-side of~\eqref{eq.nor_dyn} is continuous in~$\reals{n}$ and a solution to~\eqref{eq.nor_dyn} exists on a finite time interval~\cite[Thm. 1.1]{hal09}. Furthermore, the solutions have a continuation to a maximal interval of existence~\cite[Thm. 2.1]{hal09} and the existence of a Lyapunov function that certifies global asymptotic stability of~\eqref{eq.nor_dyn} implies that this interval is $[0,\infty)$~\cite[Prop 2.1]{bhaber00}. Such Lyapunov function is given in the proofs of Theorems~\ref{theorem.fin}~and~\ref{theorem.fxd} under the assumption that the origin is globally exponentially stable equilibrium point of~\eqref{eq.dyn}. Uniqueness of solutions follows from  the time-parameterization provided in the proof of~\cite[Prop. 1]{garbargupben22}. \end{proof}

\vsp

To recap, the existing results discussed in Section~\ref{sec.intro} establish finite/fixed-time stability for certain classes of optimization algorithms and the proofs exploit algorithm-specific inequalities. In contrast, we demonstrate that (i) any exponentially stable algorithm can be modified to guarantee convergence on a finite time interval; and (ii) finite/fixed-time stability of the modified dynamics can be certified using the Lyapunov function that proves exponential stability of the original system.




	\vspace*{-1ex}
\section{Application to optimization}
\label{sec.app}	

Next, we demonstrate the merits and the utility of the proposed approach in designing finite/fixed-time stable dynamical systems from exponentially stable optimization algorithms. We first examine nonsmooth composite optimization problems and then study smooth problems with linear constraints. To our knowledge, finite/fixed-time stability of the associated algorithms has not been previously established for either class.

	\vspace*{-1ex}
\subsection{Composite problems and related flows}\label{sec.pal}	

We consider composite optimization problems,
\beq\label{eq.comp}
\minimize\limits_x ~f(x) \,+\, g(Tx)
\eeq
where $f$: $\reals{n}\to\reals{}$ is a convex function with Lipschitz continuous gradient $\nabla f$, $g$: $\reals{d}\to\reals{}$ is a lower semicontinuous convex function, and $T\in\reals{d\times n}$ is a given matrix. Since $g$ is allowed to be nondifferentiable, constrained convex problems can be brought into the form~\eqref{eq.comp} by selecting $g$ to be the indicator function of the constraint set; see~\cite{dhikhojovTAC19} for details. 

\vsp

A solution to~\eqref{eq.comp} can be found by computing the saddle points of the associated proximal augmented Lagrangian~\cite{dhikhojovTAC19},
\beq\label{eq.pal}
\cLm(x;y)  \,\DefinedAs\, 
f(x) \;+\; \cM_{\mu g}(Tx \;+\; \mu y) \;-\; \tfrac{\mu}{2}\norm{y}_2^2.
\eeq
where $y\in\reals{d}$ is a vector of Lagrange multipliers, $\mu$ is a positive penalty parameter, and $\cM_{\mu g}$ is the Moreau envelope associated with $g$, i.e.,
\beq\non
\ba{rcl}
\!\!\!\cM_{\mu g}(v) 
	&\asp{=}& \minimize\limits_{w}~g(w) \,+\, \frac{1}{2\mu}\norm{w \,-\, v}^2
	\\[\nvs] 
&\asp{=}& g(\prox_{\mu g}(v)) \,+\, \frac{1}{2\mu}\norm{\prox_{\mu g}(v) \, - \, v}_2^2.
\ea 
\eeq
Here, the proximal operator of $g$ for a positive penalty parameter $\mu$ is defined as\cite{parboy14} 
\beq\non
\prox_{\mu g}(v)  \,=\, \argmin\limits_{w}
 g(w) \,+\, \tfrac{1}{2\mu}\norm{w \,-\, v}_2^2.
\eeq
Even for a non-differentiable $g$, the Moreau envelope is a continuously differentiable function with gradient~\cite{parboy14} 
\beq\non
\nabla\cM_{\mu g}(v) \,=\, \tfrac{1}{\mu}(v \, - \, \prox_{\mu g}(v)). 
\eeq 
In contrast to the augmented Lagrangian associated with~\eqref{eq.comp}, the proximal augmented Lagrangian~\eqref{eq.pal} is a continuously differentiable saddle function~\cite{dhikhojovTAC19}. Thus, primal-dual gradient flow dynamics can be used to compute its saddle points~\cite{ozajovCDC22},
\beq\label{eq.pal_dyn}
\ba{rcl}
\dot{z} &\asp{=}&  F(z)   \,=\,  \tbo{-\nabla_x\cLm(x;y)}{\nabla_y\cLm(x;y)} 
\\[0.4cm]
&\asp{=}& 
\tbo{-\nabla f(x)   \,-\,  T^T\nabla \cM_{\mu g}(Tx  \,+\,  \mu y)}{\mu(\nabla\cM_{\mu g}(Tx  \,+\,  \mu y)  \,-\,  y)}
\ea
\eeq
with $z \DefinedAs [x^T~y^T]^T$. 
\vsp

\begin{remark}
Unconstrained gradient, proximal gradient, and Douglas-Rachford splitting flows represent special cases of~\eqref{eq.pal_dyn}; see~\cite{mogjovAUT21} for details. 
\end{remark}

Although scaled dynamics~\eqref{eq.nor_dyn} associated with system~\eqref{eq.pal_dyn} has a convoluted form, its fixed-time stability follows from exponential stability of~\eqref{eq.pal_dyn}. In fact, Corollary~\ref{cor.pal_fxd} is a direct consequence of the analysis presented in Section~\ref{sec.conc}. 
	
\vsp

\begin{corollary}\label{cor.pal_fxd}
	Let $f$ in~\eqref{eq.comp} be a strongly convex function with a Lipschitz continuous gradient and let $T$ be a matrix with a full row rank. For any $\mu>0$, the unique saddle point of the proximal augmented Lagrangian~\eqref{eq.pal} is globally fixed-time stable equilibrium point of system~\eqref{eq.nor_dyn} resulting from the vector field $F$ in~\eqref{eq.pal_dyn} with the scaling factor~\eqref{eq.theorem.fxd_dyn.sig}.
\end{corollary}

\vsp

\begin{proof}
	The vector field $F$ in~\eqref{eq.pal_dyn} is Lipschitz continuous~\cite[Thm.~1]{dhikhojovTAC19} and system~\eqref{eq.pal_dyn} is globally exponential stable~\cite[Thm.~6]{dinjovACC19}. Since the quadratic Lyapunov function used in the proof of~\cite[Thm.~6]{dinjovACC19} has a Lipschitz continuous gradient, the condition~\eqref{eq.grad_lower} is satisfied. Thus, global fixed-time stability of system~\eqref{eq.nor_dyn} associated with~\eqref{eq.pal_dyn} follows from Theorem~\ref{theorem.fxd}.
\end{proof}

 
	\vspace*{-1ex}
\subsection{Linearly constrained problems and related flows}
We now consider convex problems with linear constraints,
\beq\label{eq.ineq}
\ba{rcl}
\minimize\limits_{x} &\asp{}& f(x)
\\[\nvs]
\subjectto &\asp{}& Ax  \,-\,  b  \,\leq\,  0
\ea
\eeq
where $f$: $\reals{n}\to\reals{}$ is a convex function with Lipschitz continuous gradient $\nabla f$, whereas $A\in\reals{d\times n}$ and $b \in \reals{d}$ are problem data. A solution to problem~\eqref{eq.ineq} is determined by the saddle points of the generalized Lagrangian~\cite{quli18, roc70a}
\beq\label{eq.gen_lag}
\cGm(x;y)  \,=\, f(x) \,+\, \sum\limits_{i \,=\, 1}^d h_\mu(a_i^Tx \,-\, b_i,y_i).
\eeq 
Here, $a_i^T$ is the $i$th row of $A$, $b_i$ is the $i$th element of the vector $b$, and $h_\mu(\cdot, \cdot)$ is the penalty function given by
\beq
h_\mu(v,w)  \,\DefinedAs\, \left\{
\ba{rl}
vw  \,+\, \mu v^2/2 & \mu v  \,+\, w   \,\geq\,  0
\\[\nvs]
-w^2/(2\mu) & \mu v  \,+\, w   \,\leq\,  0.
\ea
\right.
\eeq
Since $\cG_\mu$ is a continuously differentiable function that is convex in $x$ and concave in $y$~\cite{quli18,roc70a}, its saddle points can be computed using primal-dual gradient flow dynamics,
\beq\label{eq.pal_gen_dyn}
\ba{rcl}
\!\!\!
\dot{z} &\asp{=}&  F(z)   \,=\,  \tbo{-\nabla_x\cGm(x;y)}{\nabla_y\cGm(x;y)} 
\\[0.4cm]
\!\!\!
&\asp{=}& 
\tbo{\!\!\!-\nabla f(x)   - \displaystyle{\sum_{i \,=\, 1}^d} \max(\mu(a_i^Tx - b_i) + y_i,0)a_i\!\!\!}{(1/\mu) \displaystyle{\sum_{i \,=\, 1}^d} \max(\mu(a_i^Tx - b_i),-y_i) \rme_i}
\ea
\eeq
where $z \DefinedAs [x^T~y^T]^T$ and $\rme_i$ is the $i$th unit vector in the canonical basis of $\reals{d}$. 

\vsp


Corollary~\ref{cor.gen_lag_fxd} exploits the existing exponential stability results~\cite{quli18} for the original dynamics~\eqref{eq.pal_gen_dyn} to prove fixed-time stability of the modified system resulting from~\eqref{eq.pal_gen_dyn}.


\vsp

\begin{corollary}\label{cor.gen_lag_fxd}
	 Let $f$ in~\eqref{eq.ineq} be a strongly convex function with a Lipschitz continuous gradient and let $A$ be a matrix with a full row rank. For any $\mu>0$, the unique saddle point of the generalized Lagrangian~\eqref{eq.gen_lag} is globally fixed-time stable equilibrium point of system~\eqref{eq.nor_dyn} resulting from the vector field $F$ in~\eqref{eq.pal_gen_dyn} with the scaling factor~\eqref{eq.theorem.fxd_dyn.sig}.
\end{corollary}

\vsp

\begin{proof}
	The vector field $F$ in~\eqref{eq.pal_gen_dyn} is Lipschitz continuous~\cite[Eq.~(16)]{quli18} and a quadratic Lyapunov function certifies global exponential stability of~\eqref{eq.pal_gen_dyn}~\cite[Thm.~2]{quli18}. The rest of the proof follows identical arguments as in Corollary~\ref{cor.pal_fxd}.
\end{proof}

\vsp

\begin{remark}
When, in addition to the linear inequality constraint, there is a constraint $Cx = d$ in problem~\eqref{eq.ineq}, Corollary~\ref{cor.gen_lag_fxd} still holds if $[A^T~C^T]^T$ is a full row rank matrix; see Section~\ref{sec.ex2} and~\cite{quli18} for details. 
\end{remark}

\vsp

We note that positive orthant is a positively invariant set for primal-dual gradient flow dynamics resulting from the generalized Lagrangian~\eqref{eq.gen_lag}. Thus, if $y(0) \geq 0$, the dual variables $y (t)$ in~\eqref{eq.pal_gen_dyn} stay positive for all times~\cite{quli18}. Moreover, for nonlinear inequality constraints, the generalized Lagrangian~\eqref{eq.gen_lag} and associated dynamics~\eqref{eq.pal_gen_dyn} are still applicable~\cite{tanquli20} and, in contrast to~\eqref{eq.pal_dyn}, computations of orthogonal projections to constraint sets are not required. 

	\vspace*{-1ex}
\section{Computational experiments}\label{sec.num}

\begin{figure*}[t]
	\centering
	\begin{tabular}{r@{\hspace{-0.45 cm}}c@{\hspace{0.3 cm}} r@{\hspace{-0.45 cm}}c}
		&\subfigure[]{\label{fig.1}}
		&&
		\subfigure[]{\label{fig.2}}
		\\[-.20cm]
		\begin{tabular}{c}
			\vspace{.3cm}
			\normalsize{\rotatebox{90}{$\norm{x(t) \,-\, \xs}/\norm{\xs}$}}
		\end{tabular}
		&
		\begin{tabular}{c}
			\includegraphics[width=0.40\textwidth]{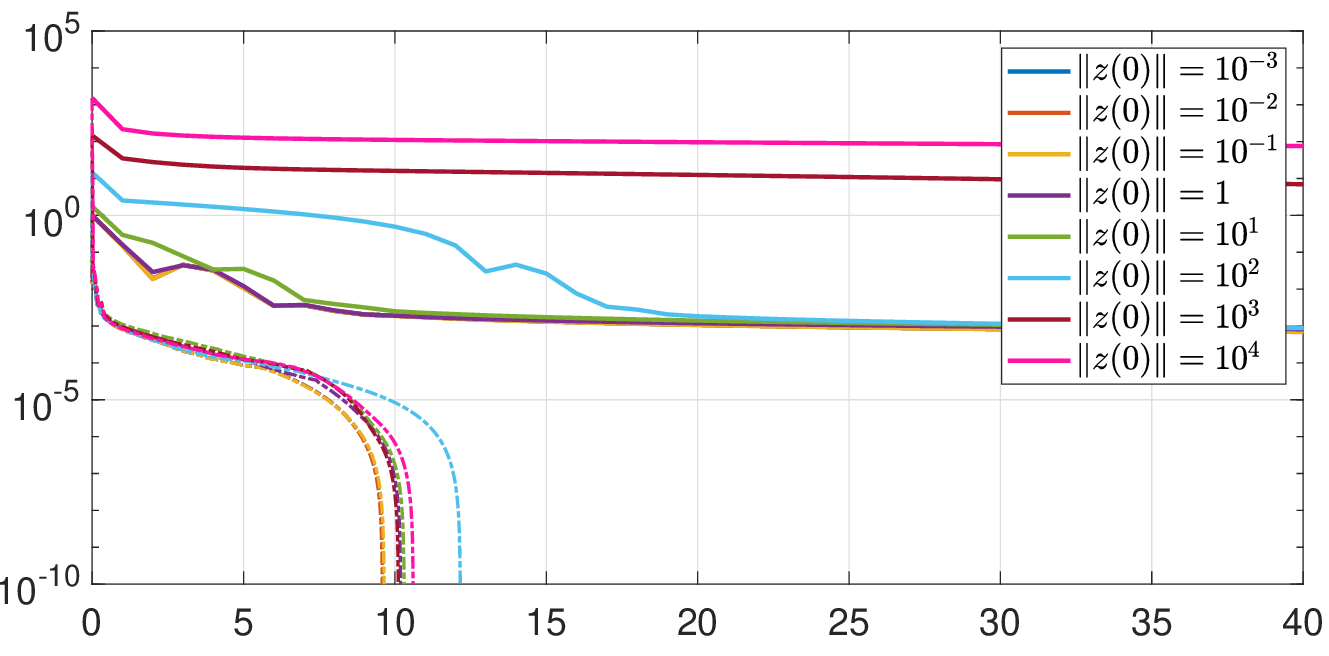}
			\\[-0.2 cm] { $t$}
		\end{tabular}
		&
		\begin{tabular}{c}
			\vspace{.3cm}
			\normalsize{\rotatebox{90}{\small $\norm{x(t) \,-\, \xs}/\norm{\xs}$}}
		\end{tabular}
		&
		\begin{tabular}{c}
			\includegraphics[width=0.40\textwidth]{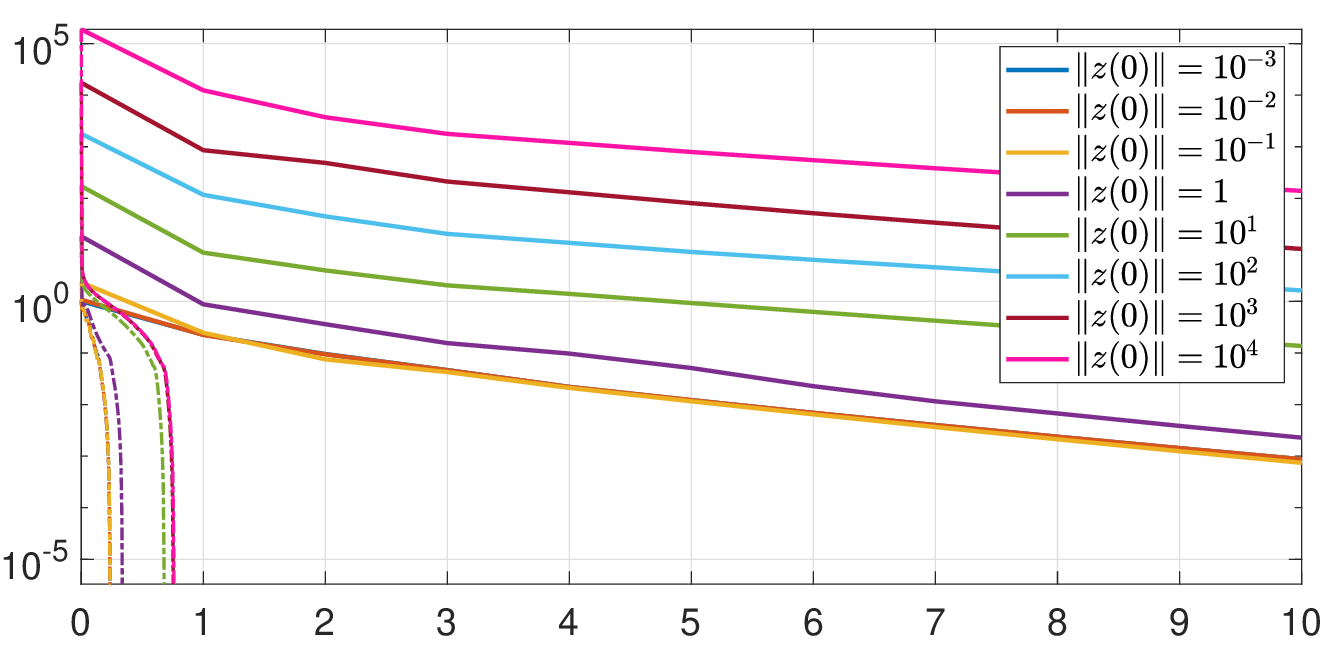}
			\\[-0.2 cm]  {$t$}
		\end{tabular}	
	\end{tabular}
	\vspace{0cm}
	\caption{Trajectories of $\zd = F(z)$ (solid lines) and $\zd = \sigma(z)F(z)$ (dashed lines) with scaling factor $\sigma$ given by~\eqref{eq.theorem.fxd_dyn.sig}. In both cases, the parameters are set to $(\mu, \eta_1, \eta_2, \lambda_1, \lambda_2) = (1, 10, 1, 1/2, 3)$ and the optimal solution $\xs$ is obtained using CVX. Different colors mark trajectories for different initial conditions for systems resulting from~(a) the proximal augmented Lagrangian~\eqref{eq.pal_dyn} associated with the total variation denoising problem~\eqref{ex.1}; and (b) the generalized Lagrangian~\eqref{eq.pal_gen_dyn} associated with the constrained quadratic program~\eqref{ex.2}.} 
	\label{fig.3}
\end{figure*}

 
We provide two examples to illustrate fixed-time stability of the modified primal-dual gradient flow dynamics resulting from proximal augmented~\eqref{eq.pal} and generalized~\eqref{eq.gen_lag} Lagrangians. Figure~\ref{fig.3} shows that the settling time of the modified dynamics~\eqref{eq.nor_dyn} with scaling factor~\eqref{eq.theorem.fxd_dyn.sig} is uniformly upper bounded irrespective of the magnitude of the initial condition.


\subsubsection{Total variation denoising}
We consider  the fused lasso~problem
\beq\label{ex.1}
\ba{rcl}
\minimize\limits_{x} &\asp{}& \frac{1}{2}\norm{Ex \,-\, q}^2 \,+\, \norm{Tx}_1
\ea
\eeq
where $E,T\in\reals{n\times n}$ and $T$ is a bidiagonal matrix with $1$ on the main diagonal and $-1$ on the first upper subdiagonal. For $n=100$, entries of $E$ are independently sampled from~$\cN(0,1)$. The measurement vector is given by $q = E \bar{x} + w$, where $w$ is noise with entries independently sampled from~$\cN(0,0.1)$ and $\bar{x}$ is a vector of random integers whose entries are independently drawn from $\set{1,\ldots,10}$ and then repeated $10$ times to get a piecewise constant partition of $\bar{x} \in \reals{n}$. 


\subsubsection{Constrained quadratic program} \label{sec.ex2}
We examine quadratic problem with linear constraints
\beq\label{ex.2}
\ba{rl}
\minimize\limits_{x}& x^TQx \,+\, q^Tx
\\[\nvs]
\subjectto&\!Ex \,=\, 0, \quad Fx \,\leq\, 0
\ea
\eeq
where $Q = R^TR$, $R\in\reals{n\times n}$,  and $E,F\in\reals{d\times n}$. We set $(d,n) = (100, 20)$ and draw entries of the matrices $R,E,F$ and the vector $q$ independently from $\cN(0, 1)$. 


	\vspace*{-1ex}
\section{Conclusion}
	\label{sec.conc}

We demonstrate that exponentially stable optimization algorithms provide a foundation for systematic design of techniques that compute the solution in finite-time. The modified algorithm is obtained by a straightforward scaling of the generator of the original dynamics and desired stability properties are certified using the Lyapunov function that proves exponential stability. Our framework departs from the analysis of specific problem instances and it can be used to achieve finite/fixed-time stability by exploiting the existing exponential stability results for a host of relevant applications including primal-dual gradient flow dynamics~\cite{quli18,dhikhojovTAC19}, anytime feasible solvers for constrained optimization~\cite{allcor23, muejor22}, and optimization-based feedback controllers~\cite{biacorpovdal21}.

	\vspace*{-1ex}
\bibliographystyle{IEEEtran}

\end{document}